\documentclass[11pt,dvipsnames]{article}
\usepackage[T1]{fontenc}
\usepackage[a4paper, margin=2.6cm, top=3cm, bottom=2.7cm]{geometry}
\usepackage{amsmath}
\usepackage{array}
\usepackage{booktabs}
\usepackage{caption}
\usepackage{fancyvrb}
\usepackage{graphicx}
\usepackage{hyperref}
\usepackage{url}
\hypersetup{breaklinks=true, colorlinks,
  urlcolor=RoyalBlue, citecolor=PineGreen, linkcolor=RoyalBlue}
\usepackage{libertine}
\usepackage{longtable}
\usepackage{newtxmath}
\usepackage{tabularx}
\usepackage{zi4}

\usepackage[english]{babel}
\usepackage{enumitem}
\usepackage{graphicx,subfigure,bm}
\usepackage[ruled,algo2e,linesnumbered]{algorithm2e}
\DontPrintSemicolon
\usepackage[utf8]{inputenc}

\usepackage[autoload=false, theme=default-plain]{jlcode}
\lstdefinestyle{jlcodeblockstyle}{%
basicstyle={\loadthemecolors\color{jlstrnum}\jlbasicfont},
keywordstyle={[1]\color{jlkeyword}\bfseries},
keywordstyle={[2]\color{jlliteral}},
keywordstyle={[3]\color{jlbuiltin}},
keywordstyle={[4]\color{jlmacros}},
keywordstyle={[5]\color{jlfunctions}},
commentstyle={\color{jlcomment}},
stringstyle={\color{jlstrnum}},
identifierstyle={\color{jlbase}},
showstringspaces=false,
upquote=true,
tabsize=4,
}
\lstdefinestyle{jlcodeboxnostyle}{%
columns=fixed,
basewidth=\bfem,
linewidth=\columnwidth,
}
\addlitjlmacros{@properties}{@properties}{11}

\usepackage{makecell}
\usepackage{listings}
\usepackage{cleveref}

\usepackage{xcolor}
\usepackage{tcolorbox}
\usepackage{listings}
\tcbuselibrary{listings}
\tcbuselibrary{xparse} 

\definecolor{light_gray}{gray}{0.92}
\definecolor{color_in}{rgb}{0.117,0.403,0.713}
\definecolor{color_out}{rgb}{0.792,0.239,0.239}

\newlength\inwd
\newcounter{ipythcntr}
\renewcommand{\theipythcntr}{\texttt{[\arabic{ipythcntr}]}}

\NewTCBListing{jinput}{O{} O{}}{%
  breakable,
  sharp corners,
  enlarge left by=\inwd,
  width=\linewidth-\inwd,
  enhanced,
  boxrule=0pt,
  colback=light_gray,
  listing only,
  top=0pt,
  bottom=0pt,
  listing engine=listings,
  listing options={
    aboveskip=1pt,
    belowskip=1pt,
    basicstyle=\footnotesize\ttfamily,
    language=julia,
    style=jlcodestyle,
    showstringspaces=false
  },
  overlay={
    \node[
      anchor=north east,
      text width=\inwd,
      font=\footnotesize\ttfamily\color{color_in},
      inner ysep=2mm,
      inner xsep=0pt,
      outer sep=0pt
      ]
      at (frame.north west)
      {\refstepcounter{ipythcntr}\ifstrempty{#1}{}{\label{#1}}\space\space In\theipythcntr:};
  }
  #2
}
\NewTCBListing{joutput}{O{} O{}}{%
  breakable,
  sharp corners,
  enlarge left by=\inwd,
  width=\linewidth-\inwd,
  enhanced,
  boxrule=0pt,
  colback=white,
  listing only,
  top=0pt,
  bottom=0pt,
  listing engine=listings,
  listing options={
    aboveskip=1pt,
    belowskip=1pt,
    basicstyle=\footnotesize\ttfamily,
    language=julia,
    style=jlcodestyle,
    showstringspaces=false
  },
  overlay={
    \node[
      anchor=north east,
      text width=\inwd,
      font=\footnotesize\ttfamily\color{color_out},
      inner ysep=2mm,
      inner xsep=0pt,
      outer sep=0pt
      ]
      at (frame.north west)
      {Out\ifstrempty{#1}{\theipythcntr}{#1}:};
  }
  #2
}
\NewTCBListing{jprint}{O{}}{%
  breakable,
  sharp corners,
  enlarge left by=\inwd,
  width=\linewidth-\inwd,
  enhanced,
  boxrule=0pt,
  colback=white,
  listing only,
  top=0pt,
  bottom=0pt,
  listing engine=listings,
  listing options={
    aboveskip=1pt,
    belowskip=1pt,
    basicstyle=\footnotesize\ttfamily,
    language=julia,
    style=jlcodestyle,
    showstringspaces=false
  },
  overlay={
    \node[
      anchor=north east,
      text width=\inwd,
      font=\footnotesize\ttfamily\color{color_out},
      inner ysep=2mm,
      inner xsep=0pt,
      outer sep=0pt
      ]
      at (frame.north west)
      {};
  }
  #1
}

\usepackage{tcolorbox}
\tcbuselibrary{skins}
\tcbuselibrary{breakable}
\tcbset{shield externalize}

\newcommand{\thetitle}{TypedMatrices.jl: An Extensible and Type-Based Matrix Collection for Julia}
\newcommand{\thetitlenote}{%
  Version of \today.
}

\newcommand{\theabstract}{%
\texttt{TypedMatrices.jl} is a Julia package to organize test matrices. By default, the package comes with a number of built-in matrices and interfaces to help users select test cases based on their properties. The package is designed to be extensible, allowing users to define their own matrix types. We discuss the design and implementation of the package and demonstrate its usage with a number of examples.}

\newcommand{\thekeywords}{%
Test matrices,
matrix collection,
matrix algorithms,
Julia%
}

\newcommand{\theauthori}{Anzhi Zhang}
\newcommand{\theorcidi}{0000-0002-2057-4373}
\newcommand{\theaffiliationi}{%
  \department{School of Computer Science}
  \institution{University of Leeds}
  \streetaddress{Woodhouse Lane}
  \city{Leeds}
  \postcode{LS2 9JT}
  \country{UK}}
\newcommand{\theemaili}{fy22az@leeds.ac.uk}

\newcommand{\theauthorii}{Massimiliano Fasi}
\newcommand{\theorcidii}{0000-0002-6015-391X}
\newcommand{\theaffiliationii}{%
  \department{School of Computer Science}
  \institution{University of Leeds}
  \streetaddress{Woodhouse Lane}
  \city{Leeds}
  \postcode{LS2 9JT}
  \country{UK}}
\newcommand{\theemailii}{m.fasi@leeds.ac.uk}

\title{\huge\bfseries\thetitle\footnote{\thetitlenote}}

\newcommand{\department}[1]{#1,}
\newcommand{\institution}[1]{#1,}
\newcommand{\streetaddress}[1]{#1,}
\newcommand{\city}[1]{#1}
\newcommand{\postcode}[1]{#1,}
\newcommand{\country}[1]{#1}
\author{\theauthori\thanks{\theaffiliationi{} (Email: \href{mailto:\theemaili}{\theemaili}, ORCID: \href{https://orcid.org/\theorcidi}{\theorcidi}).} \and%
  \theauthorii\thanks{\theaffiliationii{} (Email: \href{mailto:\theemailii}{\theemailii}, ORCID: \href{https://orcid.org/\theorcidii}{\theorcidii}).}%
}
\date{\vspace{-30pt}}

\newcommand{\R}{\ensuremath{\mathbb{R}}}
\newcommand{\Rmn}{\ensuremath{\R^{m \times n}}}
\newcommand{\Rnn}{\ensuremath{\R^{n \times n}}}

\newcommand{\matlab}{MATLAB}
\newcommand{\matrixdepot}{\texttt{MatrixDepot.jl}}
\newcommand{\typedmatrices}{\texttt{TypedMatrices.jl}}

\begin{document}

\maketitle

\vskip 40pt

\begin{abstract}
  \noindent\theabstract
  \bigskip

  \noindent\textbf{Key words.} \thekeywords{}.
\end{abstract}

\vskip 40pt

\section{Introduction}

In order to understand the accuracy of numerical algorithms for linear algebra problems, it is convenient to run these algorithms on a selection of suitably chosen test cases.
Ideally, these test problems should satisfy three requirements.
First, the they should be representative of the class of applications for which the new algorithms are designed.
Second, their solution should be known, either in closed form or to in a way that guarantees high accuracy.
Finally, they should be parametrized so that it is possible to make them arbitrarily large and/or arbitrarily difficult.

The need for good test matrices, in particular, was recognized early on by the numerical linear algebra community, and a large body of literature has been produced in the last 70 years.
To the best of our knowledge, the need to carry out numerical experiments was first highlighted by Newman and Todd~\cite{neto58}, who remark that it would be impractical to carry out rigorous error analysis of numerical algorithms for all computational problems.

Newman and Todd focus on matrix inversion, but other authors considered other linear algebra computations.
Clement~\cite{clem59}, for example, proposes a class of tridiagonal matrices for testing purposes, giving closed formulas for their inverses, determinants, and eigenvalues.
Many family of matrices expressly designed for testing purposes were proposed in the early sixties, often in the form of short communication building on each other.
A representative example is that of the Pei matrix: Pei~\cite{pei62} proposed it in 1962, and within one year LaSor~\cite{laso63} found an explicit expression for its eigenvalues and Newbery~\cite{newb63} gave an explicit method to compute its eigenvectors.
All three communications would easily fit in one page.

A summary of the test matrices proposed in the early years of digital computing can be found in the monograph by Gegory and Karney~\cite{grka69}, who collect a large number of test matrices from the literature.
Notably, their work includes the matrices gathered by Westlake in a 22-page appendix to a monograph entirely dedicated to matrix inversion~\cite{west68}.
Gregory and Karney explicitly state that the purpose of their collection is to provide numerical examples with known solutions that can be used to test algorithms for matrix inversion and computation of eigenvalues and eigenvectors.

Test matrices have considerably grown in importance since the late sixties.
Today, these test problems are used not only to gauge the accuracy of numerical algorithms, but also to benchmark their performance, to ensure that research is easily reproducible, and to conduct theoretical investigations.
Therefore, it is not surprising that several software packages have been developed to streamline the generation of test matrices and to guarantee the reproducibility of numerical results, which is of crucial importance in the mathematical sciences~\cite{dost15}.
Existing software collections and packages are surveyed in \cref{sec:related-work}.

We introduce \typedmatrices{}, a Julia package to manage large collections of test matrices by leveraging Julia's type system.
The package is released under the MIT license, is freely available on GitHub,\footnote{\url{https://github.com/TypedMatrices/TypedMatrices.jl}} and can be installed using Julia's package manager.
The collection includes 63 built-in matrices and provides a number of testing interfaces to help users generate test matrices that satisfy a specified subset of the 37 built-in properties.
The brief tutorial in \cref{sec:taste-typedmatrices} illustrates the main features and use cases.

The package is entirely implemented in Julia\footnote{\url{https://julialang.org}}~\cite{beks17}, a high-level, high-performance programming language for numerical computing.
The Julia community is highly active, and the package welcomes contributions from the community.
\typedmatrices{} builds on the features of Julia discussed in \cref{sec:background}.

\typedmatrices{} is designed to overcome the limitations of existing tools and provide a flexible environment for dynamic matrix generation.
The package is designed to be extensible and to allow users to define their own matrix types.
The main features of the package are the following.

\begin{itemize}
\item Matrices are generated dynamically, enabling substantial memory savings, especially for large-scale computations.

\item The collection leverages Julia's type system, and matrices can be stored to optimize memory usage and runtime performance.
  This design facilitates integration with Julia's multiple dispatch mechanism and provides a modular and adaptable interface.

\item The matrix entries can be of any scalar type, including user-defined types.

\item Users can easily extend the library by defining custom matrix types or generators.

\item The package supports the creation of user-defined groups, allowing users to organize matrices into logical collections.
  This feature allows researchers to save and share predefined sets of test matrices, thereby simplifying the reproducibility of numerical experiments.
\end{itemize}

The design and implementation of \typedmatrices{} are discussed in \cref{sec:features-interfaces-design-implementation}. This section shows how users can interact with the package and customize it, describes the underlying design principles, and explains the role of the Julia type system in supporting dynamic generation.

\Cref{sec:advanced-use-cases} provides some advanced use cases and examples, illustrating the collection's applicability in real-world scenarios and algorithm development.
\Cref{sec:performance} discusses the package's performance.
\Cref{sec:testing-validation} focuses on testing, validation, and performance evaluation of the collection.
\Cref{sec:conclusions} concludes with a discussion of the impact and outlines future work.

\section{Background}
\label{sec:background}

Julia~\cite{beks17} is a dynamic programming language designed to solve the two-language problems by combining the flexibility and expressiveness of high-level prototyping languages, such as Python, and the performance of lower-level languages.
Using a just-in-time (JIT) compilation, Julia delivers performance close to that of C and Fortran,\footnote{\url{https://julialang.org/benchmarks}} and Julia code can be faster than the equivalent C code in numerical linear algebra implementations~\cite{jfr22}.

A distinctive features of Julia is the use of a multiple dispatch mechanism based on parametric polymorphism.
A function is an object that maps a tuple of arguments to a return value, which can have different behaviors, and therefore different implementations, depending on what arguments it is applied to.
In Julia, functions can be overloaded, and each possible behavior of a function, which corresponds to a different implementation, is a \emph{method} of that function.
When a function is applied, Julia chooses which method to use based on the number of arguments given and on the type of all supplied arguments.
The use of multiple dispatch is a natural choice for mathematical code, and it allows programmers to write elegant, flexible, and efficient code in most scientific domains.

\section{Related Work}
\label{sec:related-work}

As far as we are aware, the first collection to be widely available in digital form was the Harwell-Boeing sparse matrix collection~\cite{dglp82,dgl89}, a set of problems for linear systems, least squares, and eigenvalue calculations from a broad range of scientific and engineering applications.
This library included earlier collections by Everstine~\cite{ever79} and George and Liu~\cite[Chapter 9]{geli81}, which were in turn based on previously developed code that was not widely available.
The collection was later re-named the Rutherford-Boeing sparse matrix collection~\cite{dgl97}, even though the latter name does not seem to have gained much popularity.

Another contemporary example, focusing again on sparse matrices, is the Matrix Market,\footnote{\url{https://math.nist.gov/MatrixMarket}}~\cite{bprb97} an online service of the US National Institute of Standards and Technology (NIST) that provides a repository for test data to be used in comparative studies of algorithms for numerical linear algebra. The collection includes nearly 500 sparse matrices, including the already mentioned Rutherford-Boeing set, as well as tools and services for matrix generation and retrieval.
The matrices can be downloaded in one of three software-independent formats.
The native format is the Matrix Market exchange format, which comes in two flavors: the \emph{coordinate format}, suitable to represent sparse matrices via explicit coordinates, and the \emph{array format}, suitable for representing dense matrices using column-major order.
Alternatively, users can opt for the Harwell-Boeing exchange format, which can save to files not only matrices, but also linear systems, potentially with a sparse coefficient matrix and/or with multiple right-hand sides, and an exact solution vector, if known.

An early software collection of matrices was developed for \matlab{} by Higham~\cite{high91m, high95m}. These codes later became part of the \matlab{} \texttt{gallery}, which is a built-in function since the R2006a release. The \matlab{} \texttt{gallery} is one of the earliest and most widely used collection and offers a comprehensive collection of classic test matrices, including many of the most interesting matrices discussed by Gregory and Karney~\cite{grka69}.

A more recent collection is the SuiteSparse Matrix
Collection\footnote{\url{https://sparse.tamu.edu}}~\cite{dahu11,kabd19} previously known as the University of Florida Sparse Matrix Collection, which at the time of writing includes 2,893 sparse matrices which are routinely used for testing and benchmarking sparse matrix algorithms. The collection includes the Rutherford-Boeing set, together with matrices from a broad range of applications, such as electromagnetics, semiconductors, structural engineering, computational fluid dynamics, and network science.

Some collection were explicitly developed to test specific implementations of numerical software.
An early example is the matrix collection gathered by Bai et al.~\cite{bddd97}, which includes a set of problems designed to test solvers for nonsymmetric eigenvalue problems.
More recently, Marques et al.~\cite{mvdp08} have presented a collection of matrices explicitly designed to test the symmetric tridiagonal eigensolvers implemented as part of LAPACK~\cite{abbb99}.
Test matrices have also been developed with the aim of validating runs of the the HPL benchmark\footnote{\url{https://www.netlib.org/benchmark/hpl}}~\cite{fahi21a}, which is used to rank the machines in the Top500 list,\footnote{\url{https://top500.org}} and for the HPL-MxP benchmark\footnote{\url{https://hpl-mxp.org}}~\cite{fahi21b}.

Collections of problems from specific applications are also common, with network science being a prime example.
The controllable test matrix toolbox for \matlab{} (CONTEST)\footnote{\url{https://www.maths.ed.ac.uk/~dhigham/CONTEST_package.html}}~\cite{tahi09} is a collection of functions to generate adjacency matrices of random graphs and to operate on them.
Other examples include the Stanford Large Network Dataset Collection\footnote{\url{https://snap.stanford.edu/data}}, part of the Stanford Network Analysis Platform (SNAP)~\cite{leso16}, whose 68 matrices are included in the SuiteSparse sparse matrix collection since 2010, and the Network Repository\footnote{\url{https://networkrepository.com}}~\cite{roah15}, an interactive network repository that collects the adjacency matrices of thousands of networks, available in Matrix Market format.

Other application-specific collections provide linear algebra problems that include several matrices and vectors.
For example, NLEVP\footnote{\url{https://github.com/ftisseur/nlevp}}~\cite{nlevp,hnt19} offers 52 nonlinear eigenvalue problems, each including several matrices.
Many toolboxes for inverse problems provide test problems including both a coefficient matrix and a right-hand side---examples are the Regularization Tools toolbox\footnote{\url{https://www.imm.dtu.dk/\textasciitilde pcha/Regutools}}~\cite{hans94,hans07}, the IR Tools toolbox\footnote{\url{https://github.com/jnagy1/IRtools}}~\cite{ghn19}, and the AIR Tools II toolbox\footnote{\url{https://github.com/jakobsj/AIRToolsII}}~\cite{hajo18}.

The wealth of matrices available in these software packages and the heterogeneity of sources poses two challenges to users.
The first challenge is the lack of a common interface.
The toolboxes described above all provide interfaces that are relatively easy to use, but these interfaces are very different, which makes it difficult for a single project to use matrices from more than one source.
The second challenge is due to the sheer number of test matrices available, which complicates the task of finding examples that satisfy one or more properties of interest.

To tackle these two issues, Higham and Mikaitis developed Anymatrix\footnote{\url{https://github.com/north-numerical-computing/anymatrix}}~\cite{himi21}, an extensible toolbox for \matlab{}.
Anymatrix provides a unified interface to several existing matrix collections, and users can extend the toolbox by creating new groups as git repositories.
Another important feature of Anymatrix is the ability to search matrices by property and to create sets of matrices from different built-in or user-defined collections.
Anymatrix includes 146 built-in matrices, organized in 7 groups, with 49 recognized properties, and several other groups have been made available by researchers.

In Julia, test matrix collections can be organized using \matrixdepot{}\footnote{\url{https://github.com/JuliaLinearAlgebra/MatrixDepot.jl}}~\cite{zhhi16}, a porting of the \matlab{} \texttt{gallery} function with a simple grouping mechanism to handle matrix properties.
\matrixdepot{} comes with 59 built-in matrices, 13 groups, and an interface to download matrices from the SuiteSparse sparse matrix collection as well as the Matrix Market.
To some extent, this package leverages Julia's capabilities, such as multiple dispatch, to generate matrices whose elements can have any scalar data type.
However, \matrixdepot{} does not support dynamic matrix generation, and this can significantly impact memory efficiency and computational scalability.
For example, all generated matrices are stored using the default matrix types, which implies that all non-zero entries must be stored explicitly, even when they could be easily computed on the fly using a simple algebraic formula.
The use of the default matrix types also prevents the use of multiple dispatch for matrix manipulation: for example, even when explicit formulas for the determinant, inverse, norm, condition number, or eigenvalues of a matrix are known, these quantities must still be calculated using the general algorithms, because any additional information about the specific matrix is lost upon generation.

The features of the packages described in this section are summarized in \cref{table:comparison}, which includes \typedmatrices{} for ease of comparison.

\begin{table}
  \centering
  \caption{Features of matrix collection implementations. The table does not list \emph{dynamic matrix generation} and \emph{use of type system}, as \typedmatrices{} is the only matrix collection to have these features.}
  \begin{tabular}{lllccc}
    \toprule
    \multicolumn{1}{c}{Collection}
    &\multicolumn{1}{c}{Language}
    &\multicolumn{1}{c}{Location}
    &\makecell{Search by\\Properties}
    &Extensible \\
    \midrule
    Regularization Tools      & \matlab{} & Own Website     & No                  & No           \\
    IR Tools                  & \matlab{} & Own Website     & No                  & No           \\
    AIR Tools II              & \matlab{} & Own Website     & No                  & No           \\
    CONTEST                   & \matlab{} & Own Website     & No                  & No           \\
    NLEVP                     & \matlab{} & GitHub          & Yes                 & No           \\
    Matrix Market             & Language-independent & Own Website     & Yes                 & No           \\
    SuiteSparse               & Language-independent & Own Website     & Yes                 & No           \\
    \texttt{gallery}          & \matlab{} & \matlab{}       & No                  & No           \\
    Anymatrix                 & \matlab{} & GitHub          & Yes                 & Yes          \\
    \matrixdepot{}            & Julia     &GitHub          & Yes                 & No           \\
    \typedmatrices{}          & \textbf{Julia} & \textbf{GitHub} & \textbf{Yes}        & \textbf{Yes} \\
    \bottomrule
  \end{tabular}
  \label{table:comparison}
\end{table}

\section{A Taste of \typedmatrices{}}
\label{sec:taste-typedmatrices}

\typedmatrices{} requires Julia version 1.6 or higher.
The package is registered in \emph{General}, the default Julia package registry,\footnote{\url{https://github.com/JuliaRegistries/General}} and can be installed  by running the following command in the Julia REPL.

\begin{jinput}
pkg> add TypedMatrices
\end{jinput}

A quick-start guide is included in the documentation.
The following example showcases the basic use of the package.

In order to use \typedmatrices{} , users must load it with the \jlinl{using} keyword:

\begin{jinput}
using TypedMatrices
\end{jinput}

The \jlinl{list_matrices} function lists all matrices available in the package.
This function provides a quick overview of the built-in matrices, and can be used to explore the matrices available in the  collection.
Users can rely on the online documentation to obtain more information on a specific matrix.

\begin{jinput}
list_matrices()
\end{jinput}
\begin{joutput}
63-element Vector{Type{<:AbstractMatrix}}:
 Orthog
 Binomial
 Lotkin
 Randjorth
 Dorr
 ⋮
 Triw
 ChebSpec
 Chow
 Forsythe
 JordBloc
\end{joutput}

Calling the \jlinl{Hilbert} constructor with a size argument generates a Hilbert matrix of the specified size.

\begin{jinput}
A = Hilbert(5)
\end{jinput}
\begin{joutput}
5×5 Hilbert{Rational{Int64}}:
  1    1//2  1//3  1//4  1//5
 1//2  1//3  1//4  1//5  1//6
 1//3  1//4  1//5  1//6  1//7
 1//4  1//5  1//6  1//7  1//8
 1//5  1//6  1//7  1//8  1//9
\end{joutput}

The function \jlinl{properties} can be used to list the properties of a matrix.

\begin{jinput}
properties(Hilbert)
\end{jinput}
\begin{joutput}
5-element Vector{Property}:
 Property(:symmetric)
 Property(:inverse)
 Property(:illcond)
 Property(:posdef)
 Property(:totpos)
\end{joutput}

This function also acceptes a matrix instance.

\begin{jinput}
properties(A)
\end{jinput}
\begin{joutput}
5-element Vector{Property}:
 Property(:symmetric)
 Property(:inverse)
 Property(:illcond)
 Property(:posdef)
 Property(:totpos)
\end{joutput}

To view all available properties, use the \jlinl{list_properties} function.

\begin{jinput}
list_properties()
\end{jinput}
\begin{joutput}
37-element Vector{Property}:
 Property(:illcond)
 Property(:rankdef)
 Property(:orthogonal)
 Property(:totnonneg)
 Property(:positive)
 ⋮
 Property(:involutory)
 Property(:defective)
 Property(:regprob)
 Property(:toeplitz)
 Property(:tridiagonal)
\end{joutput}

All builtin matrices are in the builtin group, and a group for user-defined matrices is also available. A list of all available groups is returned by the \jlinl{list_groups} function.

\begin{jinput}
list_groups()
\end{jinput}
\begin{joutput}
2-element Vector{Group}:
 Group(:user)
 Group(:builtin)
\end{joutput}

By supplying arguments, the \jlinl{list_matrices} function can be used to search for matrices that satisfy specific properties or groups.

\begin{jinput}
list_matrices(Group(:builtin))
\end{jinput}
\begin{joutput}
63-element Vector{Type{<:AbstractMatrix}}:
 Orthog
 Binomial
 Lotkin
 Randjorth
 Dorr
 ⋮
 Triw
 ChebSpec
 Chow
 Forsythe
 JordBloc
\end{joutput}
\begin{jinput}
list_matrices(Property(:symmetric))
\end{jinput}
\begin{joutput}
21-element Vector{Type{<:AbstractMatrix}}:
 DingDong
 Hankel
 Cauchy
 Poisson
 Randcorr
 ⋮
 Wilkinson
 Minij
 InverseHilbert
 Ipjfact
 RIS
\end{joutput}

Users can search for matrices with a more complex set of properties.
The following example demonstrates how to list all built-in matrices that satisfy the \verb+inverse+, \verb+illcond+, and \verb+eigen+ properties.

\begin{jinput}
list_matrices(
    [
        Group(:builtin),
    ],
    [
        Property(:inverse),
        Property(:illcond),
        Property(:eigen),
    ]
)
\end{jinput}
\begin{joutput}
3-element Vector{Type{<:AbstractMatrix}}:
 Lotkin
 Forsythe
 Pascal
\end{joutput}

By default, calling \jlinl{list_matrices} with symbols as arguments will search for matrices that satisfy these properties.
For example, the following command lists all symmetric positive definite matrices with known eigenvalues.

\begin{jinput}
list_matrices(:symmetric, :eigen, :posdef)
\end{jinput}
\begin{joutput}
4-element Vector{Type{<:AbstractMatrix}}:
 Minij
 Pascal
 Poisson
 Wathen
\end{joutput}

\section{Features, Interfaces, Design, and Implementation}
\label{sec:features-interfaces-design-implementation}

\typedmatrices{} leverages Julia's type system and multiple dispatch to generate matrices dynamically and support the use of specialized implementations of linear algebra routines.
The package is designed to be flexible and easily extensible, and its implemented aims to support high performance.
This section describes the design and implementation of the package, focusing on its core features and interfaces.

\subsection{Matrices}

Each class of matrices available in \typedmatrices{} corresponds to a Julia type.
For example, the \texttt{Hilbert} and \texttt{Clement} types represents Hilbert and Clement matrices, respectively.
Each type is a subtype of an abstract matrix type, which defines the common interfaces for all matrices in the collection.
This design ensures that the properties of a matrix are integrated with the underlying data structure and can therefore be used by specialized versions of the functions.

At present, the package includes 63 matrices with 37 properties (\cref{table:properties}).
Most of the matrices available are classic test matrices from the literature, such as Higham's test matrices~\cite{high91m}, Hansen's regularization problems~\cite{hans94,hans07}, and randomly generated matrices with specific properties.
Most of the matrices available in \matlab{} are also included, although some of those provided in the \matlab{} \texttt{gallery} are not well suited to the type system of \typedmatrices{} and are therefore not currently available---details can be found on a dedicated issue on GitHub.\footnote{\url{https://github.com/TypedMatrices/TypedMatrices.jl/issues/3}}

The package relies extensively on the multiple dispatch capability of Julia, and instead of providing matrix generators, as done in \matrixdepot{} or Anymatrix, it overloads the \jlinl{getindex} function to generate matrix entries on the fly when required.
Multiple dispatch is also used to specialize the functions in \jlinl{LinearAlgebra.jl}: when known, closed formulas are used to compute, for example, the inverse, the determinant, or the eigenvalues and eigenvectors of a test matrix, and properties are used to determine whether a matrix is symmetric, Hermitian, positive definite, and so on.

Every matrix is implemented as a new Julia type.
For example, the Hilbert matrix $H \in \Rmn$, whose elements are defined by
\begin{equation}
  \label{eq:hilbert}
  h_{i j}=\frac{1}{i+j-1},
\end{equation}
corresponds to a new Julia type, \jlinl{Hilbert}.
Here is the definition of the type in \typedmatrices{}, with its default constructor that only requires the two dimensions $m$ and $n$.
\begin{lstlisting}[language=julia, style=jlcodestyle]
struct Hilbert{T<:Number} <: AbstractMatrix{T}
    m::Integer
    n::Integer

    function Hilbert{T}(m::Integer, n::Integer) where {T<:Number}
        m >= 0 || throw(ArgumentError("$m < 0"))
        n >= 0 || throw(ArgumentError("$n < 0"))
        return new{T}(m, n)
    end
end
\end{lstlisting}
The matrix entries are not explicitly computed, and the formula~\eqref{eq:hilbert} is not present at all in the type definition: in fact, generating the matrix only assigns a value to the number of rows and columns.
Therefore, the time needed to generate a matrix of type \jlinl{Hilbert} and the memory necessary to store it do not depend on the size of the matrix.
This is in stark constant with \matrixdepot{}.

The following example shows the storage and computational efficiency of \typedmatrices{}, compared with \matrixdepot{}, to generate a square Hilbert matrix of order 100.

\begin{jinput}
@btime Hilbert(100) # TypedMatrices.jl
@btime matrixdepot("hilb", 100) # MatrixDepot.jl
sizeof(Hilbert(100)), sizeof(matrixdepot("hilb", 100))
\end{jinput}
\begin{joutput}
0.900 ns (0 allocations: 0 bytes)
12.400 μs (39 allocations: 80.17 KiB)
sizeof(Hilbert(100)) = 16
sizeof(matrixdepot("hilb", 100)) = 80000
\end{joutput}

Since \jlinl{Hilbert} is a subtype of \jlinl{AbstractMatrix}, users can interact with this matrix using the functions defined in \jlinl{LinearAlgebra}.
This is achieved by overloading the \jlinl{getindex} function for the \jlinl{Hilbert} type, so that the elements are generated dynamically, when needed, using~\eqref{eq:hilbert}.

\begin{lstlisting}[language=julia, style=jlcodestyle]
function getindex(A::Hilbert{T}, i::Integer, j::Integer) where {T}
    @boundscheck checkbounds(A, i, j)
    return T(one(T) / (i + j - 1))
end
\end{lstlisting}

Matrix types typically have a default constructor, which only requires the matrix dimensions.
The following example shows how to generate a $3 \times 3$ Hilbert matrix.
The default element type for this matrix is \jlinl{Rational\{Int64\}}.

\begin{jinput}
Hilbert(3)
\end{jinput}
\begin{joutput}
3×3 Hilbert{Rational{Int64}}:
  1    1//2  1//3
 1//2  1//3  1//4
 1//3  1//4  1//5
\end{joutput}

Several matrix types provide additional constructors for convenience or customization.
For example, the Companion matrix can be generated from a vector of coefficients, a polynomial, or a scalar.

\begin{jinput}
methods(Companion)
\end{jinput}
\begin{joutput}
\# 3 methods for type constructor:
 [1] Companion(v::AbstractVector{T}) where T<:Number
 [2] Companion(n::Integer)
 [3] Companion(polynomial::Polynomials.Polynomial)
\end{joutput}

Users can define the type of the entries of the matrix to be generated.
The following code demonstrates how the $3 \times 3$ Hilbert matrix with binary64 floating-point entries can be generated.

\begin{jinput}
display(Hilbert{Rational{Int64}}(3, 3))
display(Hilbert{Float64}(3, 3))
# here we see that element type is different, but the values are the same
Hilbert{Rational{Int64}}(3, 3) == Hilbert{Rational{Int64}}(3) == Hilbert(3, 3) == Hilbert(3) ≈ Hilbert{Float64}(3, 3)
\end{jinput}
\begin{joutput}
3×3 Hilbert{Rational{Int64}}:
  1    1//2  1//3
 1//2  1//3  1//4
 1//3  1//4  1//5
3×3 Hilbert{Float64}:
 1.0       0.5       0.333333
 0.5       0.333333  0.25
 0.333333  0.25      0.2
true
\end{joutput}

The elements of the inverse of the Hilbert matrix are known in closed form: if $H \in \Rnn$ is a Hilbert matrix, then
\begin{equation}
  \label{eq:invhilb}
  (H^{-1})_{ij} = (-1)^{i+j}(i+j-1)\binom{n+i-1}{n-j}\binom{n+j-1}{n-i}\binom{i+j-2}{i-1}^{2}.
\end{equation}

The package provides a type for the inverse of a Hilbert matrix, \jlinl{InverseHilbert}, and~\eqref{eq:invhilb} is used to define the \jlinl{getindex} function for this type.
For the \jlinl{Hilbert} type, the function \jlinl{inv} is overloaded and returns a matrix of type \jlinl{InverseHilbert} of appropriate size.

\begin{jinput}
inv(Hilbert(3))
\end{jinput}
\begin{joutput}
3×3 InverseHilbert{Int64}:
   9   -36    30
 -36   192  -180
  30  -180   180
\end{joutput}

Being subtypes of \jlinl{AbstractMatrix}, the new matrix can leverage Julia's linear algebra ecosystem.
By relying on the generic implementations in \jlinl{LinearAlgebra}, for example, users can directly compute the inverse, determinant, norm, or eigendecomposition of any of the matrices in \typedmatrices{}.

\begin{jinput}
A = Hilbert(3)
isdiag(A), isposdef(A), issymmetric(A), norm(A)
\end{jinput}
\begin{joutput}
(false, true, true, 1.413624183909335)
\end{joutput}

Where possible, the package optimizes these operations using known formulas or more efficient algorithms.
For instance, the determinant of a Hilbert matrix is computed using the closed form, which is a special case of the formula for the Cauchy determinant.

\begin{lstlisting}[language=julia, style=jlcodestyle]
LinearAlgebra.det(A::Hilbert) = inv(det(inv(A)))
\end{lstlisting}

\begin{jinput}
det(Hilbert(3))
\end{jinput}
\begin{joutput}
0.000462962962962963
\end{joutput}

Users can also convert the matrix to a variable of type \jlinl{Matrix} and compute the determinant of the latter.

\begin{jinput}
det(convert(Matrix{Float64}, Hilbert(3)))
\end{jinput}
\begin{joutput}
0.00046296296296296125
\end{joutput}

Another example is the Minij matrix~\cite{focu97}.
It is well known that its eigenvalues are of the form
\begin{equation*}
  \lambda_i = \frac{1}{4} \sec^2\left(\frac{i\pi}{2n+1}\right),\qquad i = 1, 2, \ldots, n.
\end{equation*}
Therefore, a new Julia \jlinl{LinearAlgebra} function can be defined to compute the eigenvalues of the \texttt{Minij} matrix.

\begin{lstlisting}[language=julia, style=jlcodestyle]
LinearAlgebra.eigvals(A::Minij) = [0.25 * sec(i * π / (2 * A.n + 1))^2 for i = 1:A.n]
\end{lstlisting}

\subsection{Properties}

Matrix properties play a central role in \typedmatrices{}, as they enable users to tag matrices and search the collection for examples that satisfy supplied mathematical or structural properties.
The 37 built-in properties currently available are listed in \cref{table:properties}.

\begin{table}
  \centering
  \begin{tabular}{ll}
    \toprule
    Property    & Description                                                                \\
    \midrule
    Bidiagonal  & The matrix is upper or lower bidiagonal.                                   \\
    Binary      & The matrix has entries from a binary set.                                  \\
    Circulant   & The matrix is circulant.                                                   \\
    Complex     & The matrix has complex entries.                                            \\
    Correlation & The matrix is a correlation matrix.                                        \\
    Defective   & The matrix is defective.                                                   \\
    DiagDom     & The matrix is diagonally dominant.                                         \\
    Eigen       & Part of the eigensystem of the matrix is explicitly known.                 \\
    FixedSize   & The matrix is only available in some fixed sizes.                          \\
    Graph       & The matrix is the adjacency matrix of a graph.                             \\
    Hankel      & The matrix is a Hankel matrix.                                             \\
    Hessenberg  & The matrix is an upper or lower Hessenberg matrix.                         \\
    IllCond     & The matrix is ill-conditioned for some parameter values.                   \\
    Indefinite  & The matrix is indefinite for some parameter values.                        \\
    InfDiv      & The matrix is infinitely divisible.                                        \\
    Integer     & The matrix has integer entries.                                            \\
    Inverse     & The inverse of the matrix is known explicitly.                             \\
    Involutory  & The matrix is involutory for some parameter values.                        \\
    Nilpotent   & The matrix is nilpotent for some parameter values.                         \\
    NonNeg      & The matrix is nonnegative for some parameter values.                       \\
    Normal      & The matrix is normal.                                                      \\
    Orthogonal  & The matrix is orthogonal for some parameter values.                        \\
    Positive    & The matrix is positive for some parameter values.                          \\
    PosDef      & The matrix is positive definite for some parameter values.                 \\
    Random      & The matrix has random entries.                                             \\
    RankDef     & The matrix is rank deficient.                                              \\
    Rectangular & The matrix is rectangular for some parameter values.                       \\
    RegProb     & The output is a test problem for regularization methods.                   \\
    SingVal     & Part of the singular values and vectors of the matrix is explicitly known. \\
    Sparse      & The matrix is sparse.                                                      \\
    Symmetric   & The matrix is symmetric for some parameter values.                         \\
    Triangular  & The matrix is upper or lower trinagular.                                   \\
    Tridiagonal & The matrix is tridiagonal.                                                 \\
    Toeplitz    & The matrix is Toeplitz.                                                    \\
    TotNonNeg   & The matrix is totally nonnegative for some parameter values.               \\
    TotPos      & The matrix is totally positive for some parameter values.                  \\
    Unimodular  & The matrix is unimodular for some parameter values.                        \\
    \bottomrule
  \end{tabular}
  \caption{Properties in \typedmatrices{}.}
  \label{table:properties}
\end{table}

In \typedmatrices{} each property is defined as new types.
The following code, for example, defines the \jlinl{Symmetric} property.

\begin{lstlisting}[language=julia, style=jlcodestyle]
struct Symmetric <: AbstractProperty end
\end{lstlisting}

The \jlinl{@properties} macro is used to assign properties to a newly defined matrix type. For instance, if the matrix \jlinl{MyMatrix} is symmetric and positive definite, one can associate these two predefined properties with this macro.

\begin{jinput}[code:properties_registration]
@properties MyMatrix [:symmetric, :posdef]
\end{jinput}

This is implemented by overloading the \jlinl{properties} function: the \jlinl{@properties} macro will use the \jlinl{@eval} macro to define a new \jlinl{properties} function that returns the properties associated with the matrix type given as argument. The snipeet above, for example, will generate the following code.

\begin{lstlisting}[language=julia, style=jlcodestyle]
properties(::MyMatrix) = [PropertyTypes.Symmetric, PropertyTypes.PositiveDefinite]
\end{lstlisting}

A user can list all the properties associated with a given matrix by using the \jlinl{properties} function, which returns a list of all the properties associated with that type.

\begin{jinput}
properties(MyMatrix)
\end{jinput}
\begin{joutput}
2-element Vector{Property}:
 Property(:symmetric)
 Property(:posdef)
\end{joutput}

When using the package interfaces, a property can be used as either a \jlinl{Symbol} (e.g., \jlinl{:symmetric}) or a type (e.g., \jlinl{PropertyTypes.Symmetric}). Internally, the package uses types to process properties, but the user can use symbols to make the code more readable---this will be discussed in later sections.

\subsection{Grouping}
\label{sec:grouping}

Users can organize logical collections and group matrices, for example, by properties or by experimental needs.
These groups are subsets of matrix types that can be modified dynamically and can be shared with other users of \typedmatrices{}, to facilitate reproducible research and collaborative workflows.

By default, the package includes two built-in groups: \texttt{builtin} and \texttt{user}.
The \texttt{builtin} group contains all test matrices that come with the package and are available to the user out of the box.
The \texttt{user} group is a recommended default group to which users can add custom-defined matrices.
In addition to these default groups, users can create new groups as needed and dynamically add or remove matrices.

The package provides a set of intuitive interfaces to manage groups.

\begin{itemize}
  \item \jlinl{list\_groups()}: Lists all available groups.
  \item \jlinl{add\_to\_groups(matrix, group)}: Adds a matrix to a specified group.
  \item \jlinl{remove\_from\_group(matrix, group)}: Removes a matrix from a specified group.
  \item \jlinl{remove\_from\_all\_groups(matrix)}: Removes a matrix from all groups.
  \item \jlinl{save\_group(group, file_name)}: Save group to file.
  \item \jlinl{load\_group(group, file_name)}: Load group from file.
\end{itemize}

Adding a matrix to a non-existing group creates that group, and removing all matrices from a group deletes that group, as groups cannot be empty.
In the example below, a matrix is added to the group \texttt{mygroup}, the group is saved to file, deleted, and its content is loaded to a new group \texttt{mynewgroup}.

\begin{jinput}[code:grouping_example]
add_to_groups(Minij, :user, :mygroup)   # Creates the group :mygroup.
display(list_matrices(Group(:mygroup)))
save_group(:mygroup, "mygroup.grp")
remove_from_group(Minij, :mygroup)      # Deletes group, now empty.
load_group(:mynewroup, "mygroup.grp")
display(list_groups())
display(list_matrices(:mynewgroup))
\end{jinput}
\begin{joutput}
1-element Vector{Type{<:AbstractMatrix}}:
 Minij
4-element Vector{Group}:
 Group(:mynewgroup)
 Group(:builtin)
 Group(:user)
1-element Vector{Type{<:AbstractMatrix}}:
 Minij
\end{joutput}

\subsection{Testing Interfaces}

\typedmatrices{} provides a testing framework designed to streamline the evaluation of numerical algorithms. Users can run a test function on all matrices of given dimensions that satisfy a given set of properties.

The \jlinl{test_algorithm} function takes an algorithm and optionally accepts a set of matrix properties, size constraints, a list of excluded matrices, and various options to handle errors.
The package automatically constructs all suitable matrices and applies the algorithm to them.
The following code computes the sum of the entries in $A^{2}$ for all symmetric positive definite $A$ of size between 1 and 4 available in the \jlinl{builtin} group of \typedmatrices{}.

\begin{jinput}
using BenchmarkTools

function my_function(A::AbstractMatrix)
    bench = @benchmark sum(A * A)
    return sum(A), median(bench.times)
end

test_algorithm(my_function, Integer[1,2,3,4],
    props=[Property(:posdef), Property(:eigen)],
    ignore_errors=true)
\end{jinput}

The code generates a report containing the matrix type, the size, and the output of \jlinl{my_function}.

\begin{joutput}
14-element Vector{Any}:
 (Wathen, 1, (297.77849642299265, 256375.0))
 (Wathen, 2, (820.9979833608267, 250042.0))
 (Wathen, 3, (1710.2298321182072, 251292.0))
 (Wathen, 4, (2579.9152224400186, 250083.0))
 (Pascal, 1, (1, 250958.0))
 (Pascal, 2, (5, 240958.0))
 (Pascal, 3, (19, 236167.0))
 (Pascal, 4, (69, 236583.0))
 (Minij, 1, (1, 235500.0))
 (Minij, 2, (5, 235875.0))
 (Minij, 3, (14, 250917.0))
 (Minij, 4, (30, 232292.0))
 (Poisson, 1, (4, 230666.0))
 (Poisson, 4, (8, 242583.0))
\end{joutput}

If a matrix of the requested size cannot be generated the function will return an error by default.
This behavior can be modified by setting to true the argument \jlinl{errors_as_warnings}, whereby errors will be converted to warnings, or the argument \jlinl{ignore_errors}, which will ignore all errors.
If both arguments are set to true, the latter takes precedence and errors are ignored without emitting a warning.

\subsection{Code Organization}

The package is divided into three main directories:
\begin{itemize}
\item \texttt{docs/}, which contains the documentation;
\item \texttt{src/}, which contains the source code; and
\item \texttt{test/}, which contains a suite of tests for validating the package functionality.
\end{itemize}

The entry point of the package is \texttt{src/}\typedmatrices{}.
The matrix definitions are located in the directory \texttt{src/matrices/}, with one file dedicated to each test matrix.
These individual, matrix-specific files are imported by \texttt{src/matrices/index.jl}, which exports all matrix types defined in the package, and declares the built-in group.

The remaining files in \texttt{src/} contain the code to manage the test matrices in the collection and their properties:
\begin{itemize}
\item \texttt{src/matrices.jl} provides functions to manage groups and search matrices by property;
\item \texttt{src/types.jl} defines the core types, including \jlinl{Property}, \jlinl{Group}, and individual types for each of the properties allowed by the package;
\item \texttt{src/metadata.jl} provides the mapping from properties to symbols, and functions to list list and register matrix properties; and
\end{itemize}

\subsection{Documentation}

The package's documentation is generated using Julia's \texttt{Documenter.jl} package\footnote{\url{https://github.com/JuliaDocs/Documenter.jl}} and is hosted on GitHub pages.\footnote{\url{https://typedmatrices.github.io/TypedMatrices.jl}}

The documentation is divided into two main sections: the manual and the API reference.
The manual provides an overview of the package, including installation instructions, usage examples, and contribution guidelines.
Notably, it features two key pages: \textit{Getting Started} and \textit{Performance}.
The \textit{Getting Started} page walks users through the installation process and demonstrates basic usage scenarios, helping new users become familiar with the package.
The \textit{Performance} page highlights the package's strengths in memory efficiency and computational speed, showcasing benchmarks and comparisons to emphasize its advantages.

The API reference offers a comprehensive list of all types, interfaces, properties, and built-in matrices included in the package. Each entry is accompanied by detailed descriptions and illustrative examples, enabling users to fully understand the functionality and potential applications of the package components.

The documentation is designed to be both informative and accessible, allowing users to quickly grasp the capabilities of \typedmatrices{} and efficiently integrate it into their workflows. By combining practical examples, performance insights, and detailed API references, the documentation serves as a valuable resource for both new and experienced users.

\section{Advanced Use Cases}
\label{sec:advanced-use-cases}

\typedmatrices{} can be integrated into more complex workflows, such as testing linear solvers, eigenvalue computations, or optimization algorithms. For instance, a researcher studying the performance of sparse solvers can leverage the package's support for sparse matrix types and properties, enabling detailed benchmarking and performance analysis.

This section demonstrates several practical advanced use cases, highlighting the package's capabilities in matrix generation, property-based searches, algorithm testing, and reproducible experiments. These examples demonstrate the flexibility and power of \typedmatrices{} in addressing diverse numerical computing challenges, and the package provides a comprehensive toolkit for advancing numerical research.

\subsection{Property-Based Searches and Algorithm Validation}
\label{sec:prop-based-search}

\typedmatrices{} enables researchers to efficiently identify and utilize matrices with specific properties.
One of the core features of the package is its ability to search for matrices based on their mathematical or structural properties.
For example, users can quickly retrieve all matrices that are both symmetric and positive definite, allowing for focused algorithm testing.
For instance, the following query lists all symmetric positive definite matrices in the collection:

\begin{jinput}
list_matrices(Property(:symmetric), Property(:posdef))
\end{jinput}
\begin{joutput}
14-element Vector{Type{<:AbstractMatrix}}:
...
\end{joutput}

Building on this, \typedmatrices{} provides an interface for testing algorithms directly on matrices that satisfy specified properties.
The \jlinl{test_algorithm} function automates this process by generating appropriate matrices, applying the user-defined algorithm, and returning detailed results.
For example, a user can define a custom algorithm that checks whether the determinant of a matrix is positive and validate this algorithm on symmetric positive definite matrices of various sizes:

\begin{jinput}
function custom_algorithm(A)
    return det(A) > 0
end

test_algorithm(custom_algorithm, props=[:symmetric, :posdef], sizes=[4])
\end{jinput}
\begin{joutput}
14-element Vector{Any}:
 (Cauchy, 4, true)
 (Poisson, 4, true)
 (Pascal, 4, true)
 (KMS, 4, true)
 (Pei, 4, true)
 (Moler, 4, true)
 (Prolate, 4, true)
 (Wathen, 4, true)
 (GCDMat, 4, true)
 (Lehmer, 4, true)
 (Hilbert, 4, true)
 (Minij, 4, true)
 (InverseHilbert, 4, true)
 (Ipjfact, 4, true)
\end{joutput}

This streamlined workflow reduces the time and effort required to prepare experiments, allowing researchers to focus on analyzing results and refining their algorithms. By integrating powerful search capabilities with flexible testing interfaces, \typedmatrices{} not only simplifies the experimentation process but also empowers researchers to address a wide range of numerical challenges. Its emphasis on reproducibility and ease of use further solidifies its role as a critical tool for advancing numerical algorithm research.

\subsection{Extending the Package with Custom Matrices}

\typedmatrices{} provides a flexible framework for creating and integrating custom matrix types.
Users can define new matrices in \typedmatrices{} by creating custom types and assigning them properties with the \jlinl{@properties} macro.
For example, the following code demonstrates how to define a custom \jlinl{Sumij} type for a matrix $A \in \Rmn$ with
\begin{equation*}
  a_{ij} = i + j.
\end{equation*}
As an example, the Sumij matrix of order 10 is illustrated in \cref{fig:sumij}.

\begin{figure}[t]
  \centering
  \begin{tikzpicture}[scale=0.7]
    \def\matrixsize{10}
    \pgfmathsetmacro{\maxvalue}{2 * \matrixsize - 2}

    \foreach \i in {0,...,\matrixsize} {
      \foreach \j in {0,...,\matrixsize} {
        \pgfmathsetmacro{\value}{\i + \j};
        \pgfmathsetmacro{\normvalue}{\value / \maxvalue};
        \pgfmathsetmacro{\rval}{min(max(4*\normvalue - 1.5, 0), 1)}
        \pgfmathsetmacro{\gval}{min(max(4*\normvalue - 0.5, 0), 1) - min(max(4*\normvalue - 2.5, 0), 1)}
        \pgfmathsetmacro{\bval}{1 - min(max(4*\normvalue - 1, 0), 1)}

        \definecolor{cellcolor}{rgb}{\rval,\gval,\bval};
        \fill[cellcolor] (\j,-\i) rectangle (\j+1,-\i-1);
      }
    }

    \foreach \i in {0,...,\matrixsize} {
      \node[xshift=0.375cm,yshift=0.25cm] at  (\i,0) {\i};
      \node[xshift=-0.25cm,yshift=-0.3755cm] at (0,-\i) {\i};
    }

  \end{tikzpicture}
  \newcommand{\mycaption}{The Minij matrix is a symmetric matrix with elements defined as $i+j$.}
  \caption{\mycaption}
  \label{fig:sumij}
\end{figure}

To add the \jlinl{Sumij} matrix to \typedmatrices{}, we create a new Julia type and implement its core functionality, including dynamic element generation and property checking.
The following code defines the \texttt{Sumij} type.

\begin{jprint}
import Base: getindex, size

struct Sumij{T<:Number} <: AbstractMatrix{T}
    m::Int
    n::Int
end

# Define type constructors.
Sumij(n::Integer) = Sumij(n, n)
Sumij(m::Integer, n::Integer) = Sumij{Int}(m, n)

# Assign known properties to type.
@properties Sumij [:symmetric, :integer, :positive]

# Overload LinearAlgebra functions.
size(A::Sumij) = (A.m, A.n)
LinearAlgebra.isdiag(A::Sumij) = A.n <= 1 ? true : false
LinearAlgebra.ishermitian(::Sumij) = true
LinearAlgebra.isposdef(::Sumij) = false
LinearAlgebra.issymmetric(::Sumij) = true
LinearAlgebra.adjoint(A::Sumij) = A
LinearAlgebra.transpose(A::Sumij) = A

# Dynamic generation of matrix elements.
function getindex(A::Sumij{T}, i::Integer, j::Integer) where {T}
    @boundscheck checkbounds(A, i, j)
    return T(i + j)
end
\end{jprint}

Once defined, the matrix can be used as any Julia matrix.
For example, the following code generates and displays $5 \times 5$ \texttt{Sumij} matrix.

\begin{jinput}[code:mymatrix_example]
Sumij(5)
\end{jinput}
\begin{joutput}
5×5 Sumij{Int64}:
 2  3  4  5   6
 3  4  5  6   7
 4  5  6  7   8
 5  6  7  8   9
 6  7  8  9  10
\end{joutput}

The new matrix can be added to the \texttt{user} as shown in~\cref{sec:grouping}.
\begin{jinput}
  add_to_groups(Sumij, :user)
\end{jinput}

\section{Performance}
\label{sec:performance}

\subsection{Linear Algebra Properties of Typed Matrices}

\texttt{LinearAlgebra.jl} provides several linear algebra operations. By utilizing the Julia type system, we can improve the performance of some of these operations for special matrices. The default method for the \jlinl{issymmetric} function, for example, checks that a matrix satisfies the definition of symmetry by accessing each matrix element. The matrix \jlinl{Minij} is known to be symmetric, and \typedmatrices{} defines a new method for \jlinl{issymmetric} that simply returns \jlinl{true}.

On the \jlinl{Minij} matrix of order 1,000, this specialized method is over 80,000 times faster than the default implementations in the median case.

\begin{jinput}
a = Minij(1000)
b = Matrix(Minij(1000))
@benchmark issymmetric(a)
@benchmark issymmetric(b)
\end{jinput}
\begin{joutput}
BenchmarkTools.Trial: 10000 samples with 999 evaluations.
Range (min … max):   9.810 ns … 89.790 ns  ┊ GC (min … max): 0.00% … 0.00%
Time  (median):     10.310 ns              ┊ GC (median):    0.00%
Time  (mean ± σ):   10.798 ns ±  2.083 ns  ┊ GC (mean ± σ):  0.00% ± 0.00%
...
BenchmarkTools.Trial: 4883 samples with 1 evaluation.
Range (min … max):  593.700 μs …  13.507 ms  ┊ GC (min … max): 0.00% … 0.00%
Time  (median):     873.400 μs               ┊ GC (median):    0.00%
Time  (mean ± σ):     1.009 ms ± 515.315 μs  ┊ GC (mean ± σ):  0.00% ± 0.00%
...
\end{joutput}

\subsection{Known Algorithm Working on \texttt{Hilbert}}

The example shows a known algorithm that works on \jlinl{Hilbert} matrices. The variable \jlinl{a} is of type \jlinl{Hilbert}, whereas \jlinl{b} is a variable of type \jlinl{Matrix} representing the same matrix. Computing the determinant of \jlinl{b} is 280 times slower and requires almost 1,000 times more memory than computing that of \jlinl{a}.

\begin{jinput}
a = Hilbert{BigFloat}(100)
b = Matrix(Hilbert{BigFloat}(100))
@benchmark det(a)
@benchmark det(b)
\end{jinput}
\begin{joutput}
BenchmarkTools.Trial: 6985 samples with 1 evaluation.
 Range (min … max):  334.500 μs … 740.291 ms  ┊ GC (min … max):  0.00% … 68.80%
 Time  (median):     564.100 μs               ┊ GC (median):     0.00%
 Time  (mean ± σ):   706.671 μs ±   8.853 ms  ┊ GC (mean ± σ):  10.32% ±  0.82%
...
BenchmarkTools.Trial: 32 samples with 1 evaluation.
 Range (min … max):  127.925 ms … 229.261 ms  ┊ GC (min … max): 8.86% … 4.76%
 Time  (median):     158.327 ms               ┊ GC (median):    9.49%
 Time  (mean ± σ):   160.932 ms ±  23.576 ms  ┊ GC (mean ± σ):  9.48% ± 3.52%
...
\end{joutput}

\subsection{Trade-off between Performance and Memory}

For algorithms that are not implemented in \typedmatrices{}, the package trades performance off for—potentially substantial—memory savings. For example, generating the variable \jlinl{a}, which is of type \jlinl{Cauchy}, only requires \textbf{63.229 $\boldsymbol{\mu}$s} and \textbf{114.16 KiB} of memory, while generating \jlinl{b}, which is the same matrix but has type \jlinl{Matrix}, requires \textbf{3.862 ms} and \textbf{7.74 MiB} of memory. And once generated, storing \jlinl{b} requires 500,000 times more memory than storing \jlinl{a}.

\begin{jinput}
@benchmark a = Cauchy{Float64}(1000)
@benchmark b = Matrix(Cauchy{Float64}(1000))
@show Base.summarysize(a)
@show Base.summarysize(b)
\end{jinput}
\begin{joutput}
BenchmarkTools.Trial: 10000 samples with 1 evaluation.
 Range (min … max):  27.100 μs … 191.819 ms  ┊ GC (min … max):  0.00% … 99.94%
 Time  (median):     32.100 μs               ┊ GC (median):     0.00%
 Time  (mean ± σ):   63.229 μs ±   1.919 ms  ┊ GC (mean ± σ):  35.05% ±  4.29%
...
BenchmarkTools.Trial: 1288 samples with 1 evaluation.
 Range (min … max):  2.413 ms … 18.386 ms  ┊ GC (min … max):  0.00% … 48.63%
 Time  (median):     3.271 ms              ┊ GC (median):     0.00%
 Time  (mean ± σ):   3.862 ms ±  1.674 ms  ┊ GC (mean ± σ):  15.96% ± 19.84%
...
Base.summarysize(a) = 16
Base.summarysize(b) = 8000040
\end{joutput}

On the other hand, accessing elements of \jlinl{a} requires computation, whereas the elements of \jlinl{b} have been pre-computed and are already available in memory. This implies a performance penalty, which is not unexpected. In view of this trade-off, however, one can use extremely large matrices on machines with a moderate amount of memory, which allows users to tackle otherwise intractably large problems. This is especially true for algorithms that only need to access a subset of the matrix elements.

\begin{jinput}
@benchmark det(a)
@benchmark det(b)
@benchmark sum(a)
@benchmark sum(b)
\end{jinput}
\begin{joutput}
BenchmarkTools.Trial: 111 samples with 1 evaluation.
 Range (min … max):  20.537 ms … 353.410 ms  ┊ GC (min … max): 0.00% … 90.93%
 Time  (median):     34.151 ms               ┊ GC (median):    0.00%
 Time  (mean ± σ):   45.104 ms ±  42.894 ms  ┊ GC (mean ± σ):  7.81% ±  9.53%
...
BenchmarkTools.Trial: 175 samples with 1 evaluation.
 Range (min … max):  18.639 ms … 314.529 ms  ┊ GC (min … max): 0.00% … 91.89%
 Time  (median):     26.317 ms               ┊ GC (median):    0.00%
 Time  (mean ± σ):   28.670 ms ±  22.610 ms  ┊ GC (mean ± σ):  7.81% ±  8.48%
...
BenchmarkTools.Trial: 3104 samples with 1 evaluation.
 Range (min … max):  1.124 ms …   7.772 ms  ┊ GC (min … max): 0.00% … 0.00%
 Time  (median):     1.400 ms               ┊ GC (median):    0.00%
 Time  (mean ± σ):   1.604 ms ± 579.750 μs  ┊ GC (mean ± σ):  0.00% ± 0.00%
...
BenchmarkTools.Trial: 10000 samples with 1 evaluation.
 Range (min … max):  243.900 μs …  2.106 ms  ┊ GC (min … max): 0.00% … 0.00%
 Time  (median):     329.800 μs              ┊ GC (median):    0.00%
 Time  (mean ± σ):   355.504 μs ± 91.684 μs  ┊ GC (mean ± σ):  0.00% ± 0.00%
...
\end{joutput}

The results also summarized in \cref{table:performance_tradeoff}.

\begin{table}[t]
  \centering
  \begin{tabular}{cccccc}
    \toprule
    Matrix (of size 1000) & Operation & Memory   & Time (median)  & Time (mean ± $\sigma$)         \\
    \midrule
    Cauchy                & det       & 16 bytes & 34.151 ms      & 45.104 ms ± 42.894 ms          \\
    Matrix(Cauchy)        & det       & 7.63 MiB & 26.317 ms      & 28.670 ms ± 22.610 ms          \\
    Cauchy                & sum       & 16 bytes & 1.400 ms       & 1.604 ms ± 579.750 $\mu$s      \\
    Matrix(Cauchy)        & sum       & 7.63 MiB & 329.800 $\mu$s & 355.504 $\mu$s ± 91.684 $\mu$s \\
    \bottomrule
  \end{tabular}
  \caption{%
    Performance trade-off between memory and time for \texttt{Cauchy} matrix.%
  }
  \label{table:performance_tradeoff}
\end{table}

\section{Testing and Validation}
\label{sec:testing-validation}

\typedmatrices{} is accompanied by a suite of tests, and continuous integration is provided via GitHub actions.
The tests encompass two levels.
On the one hand, they cover the core functionalities of the package: listing matrices and properties, managing properties and groups, and batch testing of algorithms.
On the other hand, the generation of each matrix in the package is thoroughly tested to confirm that the constructors and any method defined for the type work as expected.
\begin{enumerate}
\item For types with multiple constructors, we check that generated matrices that are expected to be equivalent are indeed so.
\item We check that the type of the generated matrices matches the type specified by the user or is the default type if none is specified.
\item We generate matrices with all possible element types and check that the these are equivalent up to rounding errors.
\item For each \typedmatrices{} type, we generate a test matrix, we convert it to an object of Julia's \jlinl{Matrix} type, and we compare the result of linear algebra operations applied to the two.
\item Matrices generated with typical parameters are compared against hard-encoded matrices to confirm the correctness of the generation logic.
\end{enumerate}

The tests are implemented using Julia \texttt{Test} package and achieve 100\% code coverage, as verified through automated testing tools.

\section{Conclusions and Future Work}
\label{sec:conclusions}

We have introduced \typedmatrices{}, a Julia package to organize collections of test matrices.
\typedmatrices{} is extensible and utilizes Julia's type system to generate the matrix entries dynamically.
It supports properties and allows users to create and share group of matrices for testing purposes.

Compared with existing alternatives, the package has several advantages.
The dynamic generation significantly reduces both memory storage matrix generation time.
The presence of a testing interface to run directly tests on matrices selected by property streamlines the process of numerical tests---to the best of our knowledge, \typedmatrices{} is currently the only package to provide this kind of functionality in any language.
Finally, the extensive use of types for matrices and properties provides a more friendly user interface that supports linting and code completion.

The only alternative currently available for Julia is \matrixdepot{}.
\matrixdepot{} is a porting of the \matlab{} \texttt{gallery} function, and therefore it does not rely on Julia's type system.
A main difference between these two Julia packages is the clear separation between properties and groups: in \matrixdepot{}, a matrix is assumed to satisfy a certain property if it belongs to a group named after that property---this mixes two concepts that should be kept logically separate.

\typedmatrices{} does not currently support static collections, such as Matrix Market~\cite{bprb97} and SuiteSparse~\cite{dahu11}.
The reason for this is that if the entries cannot be computed using a known formula, then the matrices cannot be generated dynamically and, therefore, there is no benefit in using the type system.

The general design of \typedmatrices{} can be applied to any language with a hierarchical type system.
The implementation can be simplified significantly if the language supports polymorphism: for each type, only the functions that have a special implementation will require overloading, while the others will fall back to the default implementation.
Julia is a particularly suitable language for this, as the \jlinl{LinearAlgebra} packages provides generic implementations that will work with any type that extends \jlinl{AbstractMatrix} and implements the \jlinl{getindex} function.
Any language with a similar structure for linear algebra routines is a good candidate for a porting of \typedmatrices{}.

\bibliographystyle{plain-doi}

\bibliography{references}

\end{document}